\RequirePackage{fix-cm}
\documentclass[smallextended]{svjour3}       % onecolumn (second format)
\smartqed  % flush right qed marks, e.g. at end of proof
\usepackage{savesym} \savesymbol{AND}
\usepackage{latexsym,amssymb}
\usepackage{booktabs}
\usepackage{flafter}
\usepackage{graphicx}
\usepackage{algorithm}
\usepackage{algorithmic}
\usepackage{amsmath}
\usepackage{xcolor}
\usepackage{multirow}
\usepackage{bm}
\usepackage{graphicx}
\usepackage{indentfirst}
\usepackage{etex}
\reserveinserts{18}
\usepackage{morefloats}
\begin{document}

\title{A new deflated block GCROT($m,k$) method for the solution of linear systems with
multiple right-hand sides  \thanks{Supported by YBXSZC20131070, the the National Natural Science Foundation of China (11026085, 11101071, 51175443) and the Project-sponsored by OATF (UESTC).}
}


\author{Jing Meng         \and
        Peiyong Zhu  \and
        Houbiao Li  \and
        Yanfei Jing
        %etc.
}


\institute{Jing Meng \and Peiyong Zhu  \and Houbiao Li  \and   Yanfei Jing   \at
           School of Mathematical Sciences, University of Electronic Science and Technology of China, Chengdu, Sichuan, 611731, P.R.China \\
              \email{mengmeng-erni@163.com, lihoubiao0189@163.com}           %  \\
         %  \and
         %  S. Author \at
         %     second address
}

\date{Received: date / Accepted: date}

\maketitle

\begin{abstract}
Linear systems with multiple right-hand sides arise in many applications. To solve such systems efficiently, a new deflated block GCROT($m,k$) method is explored in this paper by exploiting a modified block Arnoldi deflation. This deflation strategy has been shown to have the potential to improve the original deflation which indicates an explicit block size reduction. Incorporating this modified block Arnoldi deflation, the new algorithm can address the possible linear dependence at each iteration during the block Arnoldi procedure and avoids expensive computational operations. In addition, we analyze its main mathematical properties and prove that the deflation procedure is based on a non-increasing behavior of the singular values of the true block residual. Moreover, as a block version of GCROT($m,k$), the new method inherits the property of easy operability. Finally, some numerical examples also illustrate the effectiveness of the proposed method.
\end{abstract}

\section{Introduction}
\label{intro}
Let us consider the linear systems with $p$ right-hand sides (RHSs)
\begin{equation}\label{eq:11}
  AX=B,
\end{equation}
where $A\in \mathbb{C}^{n\times n}$ is a non-singular matrix of
large dimension, $B\in \mathbb{C}^{n\times p}$ is full rank and
$X\in \mathbb{C}^{n\times p}, (p \ll n)$. Such linear systems with
multiple RHSs arise in many applications, see, e.g., electromagnetic
scattering \cite{P.Soudais}, model reduction in circuit
simulation\cite{R. W. Freund.}, Quantum Chromo Dynamics
QCD\cite{Bloch4,Sakurai}.

To solve such systems efficiently, block Krylov subspace methods,
which are extended iterative solvers from single to multiple systems,
have been appealing. This is due to the fact that a block Krylov
subspace has a much larger search space and contains all
vectors of Krylov subspaces generated by a single linear system ($p=1$). Let $X_{0}\in \mathbb{C}^{n\times p}$ be an initial block guess and
$R_{0}=B-AX_{0}$. Then the approximate
solution $X_{k}\in \mathbb{C}^{n\times p}$ generated by
block iterative methods satisfies
\[
X_{k}-X_{0}\in \mathcal{K}_{k}(A,R_{0}),
\]
where
\begin{equation}\nonumber
\mathcal{K}_{k}(A,R_{0}) = \{\sum_{i=0}^{k-1}A^{i}R_{0}\gamma_{i}, \forall\gamma_{i}\in \mathbb{C}^{p\times p}, 0\leq i\leq k-1\}
\subset \mathbb{C}^{n\times p}
\end{equation}
is the $k$-th block Krylov subspace generated by $A$ and increasing powers of $A$ applied to $R_{0}$.
Note that each of the $p$ columns of $X_{k}$ satisfies
\begin{eqnarray}\nonumber
X_{k}(:,l)-X_{0}(:,l) & \in & \{\sum_{i=0}^{k-1}\sum_{j=1}^{p}A^{i}R_{0}(:,j)\gamma_{i}(j,l), \gamma_{i}(j,l) \in \mathbb{C}, \forall 1\leq l \leq p \} \\
&=& \sum_{j=1}^{p}\mathbb{K}_{k}(A,R_{0}(:,j)), \nonumber
\end{eqnarray}
where $\mathbb{K}_{k}(A,R_{0}(:,j))=\{R_{0}(:,j), AR_{0}(:,j), \ldots, A^{k-1}R_{0}(:,j)\} $ and $R_{0}(:,j)$ denotes the $j$-th column of $R_{0}$.
By comparison with the single right-hand side case, the solution of each linear system is sought in a richer space leading hopefully to a reduction of iteration counts. Moreover, another advantage is that block Krylov subspace algorithms are better suited to parallelism \cite{O'Leary,Li} and make better use of higher level BLAS \cite{Bai}.

The block Conjugate Gradient (BCG) is the first block iterative solver introduced by O'Leary \cite{O'Leary} and its related algorithms were proposed for parallel computers \cite{Nikishin,Simoncini}. For nonsymmetric problems, many block counterparts have been proposed, such as the block generalized minimal residual (BGMRES) method and its variants \cite{Vital,Simoncini11,Liu,M. Sadkane,R. D. da Cunha,Calandra}, the block quasi minimum residual (BQMR) method \cite{Freund2},  the block BiCGstab (BBiCGstab) algorithm\cite{Guennouni}, the block Lanczos method\cite{Guennouni1}, the block least squares (BLSQR) algorithm\cite{Karimi} and the block IDR($s$) algorithm\cite{Du}. Refer to \cite{Gutknecht} for a recent overview on block Krylov subspace methods.

In 1996, De Sturler \cite{Sturler} proposed a generalized conjugate residual method with implicit inner orthogonalization (GCRO) to address a single linear system problem.
GCRO is essentially a nested Krylov method based on the generalized conjugate residual (GCR) method, which uses GMRES as an inner method. The method allows users to select the optimal correction over arbitrary subspace. Furthermore, to minimize the error produced by discarding information, De Sturler extended this concept by providing a framework where the optimal subspace is retained from one cycle to the next. This method is called GCROT \cite{Sturler111}. However, GCROT is complicated to implement and requires five nontrivial parameters. In order to reduce the burden of determining optimal parameters, Hicken \& Zingg presented a competitive method, GCROT($m,k$) \cite{Jason}, which requires only two parameters: an inner subspace size and an outer subspace size. In this case, GCROT($m,k$) is straightforward to implement. Based on this new implementation, its block version, block GCROT($m,k$) (BGCROT($m,k$)), has been developed \cite{Jing}. It is a block method for solving multiple RHSs linear systems, but it is usual to come across the possible linear dependence of some columns of the block residuals. An attractive improvement measure consists in combining BGCROT($m,k$) with a deflation technique which is used to delete such a dependence explicitly during the block iterative procedure.

This deflation may occur at startup or in a later step \cite{Gutknecht}. To the best of our knowledge, there exist three ways to proceed.
\begin{itemize}

  \item  Initial deflation \cite{J. Langou,Gutknecht} is a simple strategy  which allows to remove the possible linearly dependent columns in the initial block residual $R_{0}$.
  \item  A second deflation strategy consists in deleting linearly and almost linearly dependent columns at each initial computation of the block residual, i.e., at the beginning of each cycle, when a restarted block Krylov subspace method is used \cite{Calandra,Gutknecht,H.Calandra}.
  \item  Deflation at each iteration deals with linearly and almost linearly dependent columns of block residuals occurring at each iteration \cite{M. Sadkane,Gutknecht,H.Calandra} in the block Arnoldi procedure.
\end{itemize}

As far as we know, the first method incorporating the third strategy is due to Robb\'e and Sadkane \cite{M. Sadkane}. Unfortunately, deflation may lead to a loss of information that slows down the convergence \cite{J. Langou}. To remedy this situation, Robb\'e and Sadkane kept the almost linearly dependent vectors and reintroduced them in the next iterations if necessary during the block Arnoldi procedure, for detail, see \cite{M. Sadkane}. We call this skill the modified block Arnoldi deflation.

Modified block Arnoldi deflation technique has shown great potential to improve the convergence and reduce computational cost for block Krylov subspace solvers \cite{M. Sadkane,H.Calandra}, in many cases, without dramatically increasing the memory requirements. Therefore, if we can combine BGCROT($m,k$) with this technique, we will have an effective method which will be able to handle the possible linearly dependent vectors in the block Krylov subspace. We call this new method deflated BGCROT($m,k$).

The main contributions of this paper can be summarized as follows. First we will derive a new deflated BGCROT($m,k$) method (DBGCROT($m,k$)) by exploiting the modified block Arnoldi deflation technique. Incorporating this deflation strategy, the new approach can address the possible linear dependence at each iteration during the block Arnoldi procedure and avoids expensive computational operations. Second, we analyze its main mathematical properties and then prove that the deflation procedure is based on a non-increasing behavior of the singular values of the block residual.
The structure of the paper is as follows. In Section 2, we describe
the deflated block GCRO method and detail the modified block Arnoldi deflation procedure.
Thereafter, its truncated block version is derived in Section 3. The
effectiveness of the proposed method is also demonstrated in Section
4. Finally, some conclusions are summarized in Section 5.

\section{A deflated block GCRO method}

In this section, our goal is to adapt the modified block Arnoldi deflation technique to improve existing block GCRO (BGCRO) method \cite{Jing,Yu}. In addition, we detail how to address linearly and almost linearly dependent vectors during the BGCRO iteration procedure.

\subsection{ The BGCRO method}

We briefly review the BGCRO method, as described in \cite{Jing,Yu}.
Firstly the block version of GCR used as an outer iteration method is given. Let
$\mathcal{U}_{k}=[U_{1},U_{2},\ldots,U_{k}]$,  $\mathcal{C}_{k}=[C_{1},C_{2},\ldots,C_{k}]\in \mathbb{C}^{n\times kp}$ be given matrices satisfying
\[
 A\mathcal{U}_{k}= \mathcal{C}_{k},
 \]
 \[
 \mathcal{C}_{k}^{H}\mathcal{C}_{k}=  I_{kp},
\]
where $U_{i}, C_{i}\in \mathcal{K}(A,R_{0}), 1\leq i\leq k$, $U_{i}$ is a block search vector and $\mathcal{K}(A,R_{0}) = \text{span}\{R_{0}, AR_{0}, A^{2}R_{0}, \ldots\}$. % and $\mathcal{I}_{k}$ is the $kp\times kp$ identity matrix.
Then we consider the following minimization problem
\[
\min_{X \in \mathcal{R}(\mathcal{U}_{k})} ||R_{0} - AX||_{F}.
\]
In order to minimize the residual over the search space $\mathcal{R}(\mathcal{U}_{k})$, the approximate solution $X_{k}$ and the corresponding residual $R_{k}\in \mathbb{C}^{n\times p}$ satisfy
\[
{X}_{k} =  X _{0} + \mathcal{U}_{k}\mathcal{C}_{k}^{H}R_{0},
 \]
\[
 {R}_{k} = R_{0} - \mathcal{C}_{k}\mathcal{C}_{k}^{H}R_{0} \ \ \text{such that} \ \  {R}_{k} \perp \mathcal{R}(\mathcal{C}_{k}).
 \]
The following question refers to how to generate ${U}_{k+1}$ and ${C}_{k+1}$ for the subsequent iteration. Ideally, we would like to choose ${U}_{k+1}= E_{k}$ with $E_{k}=X-X_{k}$. However, in general, it is not easy to get the error $E_{k}$ due to the unknown $X$. An alternative method is to choose a suitable approximation to the error $E_{k}$, which is tantamount to solving the equation
\begin{equation}\label{eq:14}
AE_{k}={R}_{k}.
\end{equation}
This is done by using an inner iteration method. In general, any block Krylov-based iterative solver (e.g., BGMRES \cite{Vital,Simoncini11,Liu,M. Sadkane,R. D. da Cunha}, BBiCGstab \cite{Guennouni}), which gives an approximate solution to
$E_{k}$, could be used as the inner method. Here, BGMRES is
considered. For preserving the orthogonality relations of GCR in the
inner algorithm (or for a faster convergence speed), De Sturler
\cite{Sturler} explored an idea of using $(I-C_{k}C_{k}^{H})A$
instead of $A$ as the operator in the Krylov method in the inner
loop for the single linear system. Following \cite{Sturler}, we also
take $(I -\mathcal{C}_{k}\mathcal{C}_{k}^{H})A$ as the concerned operator.
So solving (\ref{eq:14}) by using BGMRES (after $m$ iterations) is equivalent to solve the following minimization problem
\[
{Y}_{m} = \mathop{\arg\min}_{Y \in \mathbb{C}^{mp\times p}} ||R_{k} - (I-\mathcal{C}_{k}\mathcal{C}_{k}^{H})A\mathcal{V}_{m}Y||_{F},
\]
where $\mathcal{V}_{m}$ is an $n \times mp$ orthogonal matrix. % generated by BGMRES(m).
Since ${R}_{k} \perp \mathcal{R}(\mathcal{C}_{k})$, we get
\[
{Y}_{m} =\mathop{\min}_{Y \in \mathbb{C}^{mp\times p}}
 ||(I -\mathcal{C}_{k}\mathcal{C}_{k}^{H})(R_{k} -A\mathcal{V}_{m}Y)||_{F}.
\]
Then, in the outer loop we set
\[
U_{k+1}= (\mathcal{V}_{m}Y_{m}-\mathcal{U}_{k}\mathcal{C}_{k}^{H}A\mathcal{V}_{m}Y_{m})/||(I-\mathcal{C}_{k}\mathcal{C}_{k}^{H})A\mathcal{V}_{m}Y_{m}||_{F},
\]
\[
C_{k+1} = ((I-\mathcal{C}_{k}\mathcal{C}_{k}^{H})A\mathcal{V}_{m}Y_{m})/||(I-\mathcal{C}_{k}\mathcal{C}_{k}^{H})A\mathcal{V}_{m}Y_{m}||_{F}.
\]
To make these ideas more concrete, the pseudocode for BGCRO is presented in Algorithm 1.

\begin{algorithm}[htbp]
\caption{The BGCRO method \cite{Jing,Yu}}             %
\label{alg:Framwork}                  % 给算法一个标签，这样方便在文中对算法的引用
\begin{algorithmic}[1]
\STATE Compute $R_{0}=B-AX_{0}$\\
\FOR{$k =0,1,\ldots$}
\STATE $\sharp$ Perform $m$ steps of the BGMRES method
\STATE Compute the QR decomposition of $R_{k}$ as $R_{k} = QR$  \\
\STATE  $V_{1} =Q$, $G_{1}=[R^{H}, 0_{mp \times p}]^{H}$\\
\FOR{$j=1,\ldots, m$,}
\STATE Compute $W=AV_{j}$\\
\STATE $\sharp$ Orthogonalize $W$ against $\mathcal{C}_{k}$
\FOR {$i= 1,\ldots,k$}
\STATE $ B_{i,j}=C_{i}^{H}W$\\
\STATE  $W:=W-C_{i}B_{i,j}$\\%\mathcal{B}
\ENDFOR
\STATE $\sharp$ Orthogonalize $W$ against $\mathcal{V}_{j}$
\FOR {$i= 1,\ldots, j$}
\STATE $H_{i,j}=V_{i}^{H}W$ %\mathcal{H}
\STATE  $W:=W-V_{i}H_{i,j}$
\ENDFOR
\STATE   Compute the QR decomposition of $W$ as $W =V_{j+1}H_{j+1,j}$
\ENDFOR
\STATE Define  $\mathcal{V}_{m+1}=[V_{1},V_{2},\ldots,V_{m+1}]$
\STATE Compute $Y_{m}=\mathop{\arg\min}_{Y \in \mathbb{C}^{mp \times p} }||G_{1}- \mathcal{H}_{m}Y ||_{F}$
\STATE $\sharp$ Define new outer vectors \\
\STATE $U_{k+1}=(\mathcal{V}_{m}-\mathcal{U}_{k}\mathcal{B}_{m})Y_{m}$ with $\mathcal{B}_{m}= \mathcal{C}_{k}^{H}A\mathcal{V}_{m}$\\
\STATE $C_{k+1}= \mathcal{V}_{m+1} \mathcal{H}_{m}Y_{m}$\\
\STATE  Compute the QR decomposition of $C_{k+1}$ as $C_{k+1} = QR$
\STATE $C_{k+1}=Q, U_{k+1}=U_{k+1}/R$\\
\STATE $\sharp$ Update residual and solution\\
\STATE $R_{k+1}:=R_{k}-C_{k+1}(C_{k+1}^{H}R_{k})$\\
\STATE $X_{k+1}:=X_{k}+U_{k+1}(C_{k+1}^{H}R_{k})$
\ENDFOR
\label{code:recentEnd}
\end{algorithmic}
\end{algorithm}

In Algorithm 1, the block Gram-Schmidt procedure (from line 5 to line 19) proceeds by orthonormalizing $AV_{j}$ against $\mathcal{C}_{k}$ and $\mathcal{V}_{j}$, which constructs a block Arnoldi-like relation
\[
(I -\mathcal{C}_{k}\mathcal{C}_{k}^{H})A\mathcal{V}_{m}=\mathcal{V}_{m+1} \mathcal{H}_{m},
\]
where
\[\mathcal{H}_{m} = \left[
  \begin{array}{cccc}
    H_{1,1} & H_{1,2} & \ldots & H_{1,m} \\
    H_{2,1} & H_{2,2} & \ldots  & H_{2,m} \\
    0_{p \times p} & H_{3,2} & \ldots  & H_{3,m}\\
    \vdots & \ddots & \ddots & \vdots \\
    0_{p \times p}  &\ldots & 0_{p \times p} & H_{m+1,m} \\
  \end{array}
\right] \in \mathbb{C}^{(m+1)p \times mp}\]
is a block Hessenberg matrix.
After minimizing the Frobenius norm of the block residual (line 21 of Algorithm 1), the inner loop procedure is done.

Let $P_{\mathcal{C}_{k}} = I-\mathcal{C}_{k} \mathcal{C}_{k}^{H}$ and $A_{\mathcal{C}_{k}} = P_{\mathcal{C}_{k}} A$. In the inner loop, we consider $m$ steps of BGMRES to find the optimal approximation in the subspace $\mathcal{K}_{m}(A_{\mathcal{C}_{k}},P_{\mathcal{C}_{k}}R_{k})$. Since $\mathcal{R}(\mathcal{C}_{k}) \subset \mathcal{K}(A, R_{0})$, $(\mathcal{K}(A, R_{0}) \cap \mathcal{R}(\mathcal{C}_{k}))^{\bot}$ is an invariant subspace of $A_{\mathcal{C}_{k}}$. It allows us to compute the optimal approximation over the (global) space $ \mathcal{R}(\mathcal{U}_{k}) + A^{-1}\mathcal{K}_{m}(A_{\mathcal{C}_{k}},A_{\mathcal{C}_{k}}R_{k})$.
The following theorem summarizes the convergence properties for the BGCRO approach.% and $\mathcal{V}_{m+1}^{H}\mathcal{V}_{m+1} = I_{(m+1)p}$
\begin{theorem}
Let $A_{\mathcal{C}_{k}} = (I -\mathcal{C}_{k}\mathcal{C}_{k}^{H} )A$ and suppose that the block Krylov matrix ($R_{k}, A_{\mathcal{C}_{k}}R_{k},\ldots, A_{\mathcal{C}_{k}}^{m-1}R_{k}$) has full column rank $mp$. Using a block Arnoldi procedure in the inner loop, we obtain the following equation $A_{\mathcal{C}_{k}}\mathcal{V}_{m} =\mathcal{V}_{m+1} \mathcal{H}_{m}$. Let $Y_{m}$ be the solution of the inner BGMRES (after $m$ iterations) method:
\[
{Y}_{m} = \mathop{\arg\min}_{Y \in \mathbb{C}^{mp\times p}} ||R_{k} - A_{\mathcal{C}_{k}}\mathcal{V}_{m}Y||_{F} =
\mathop{\arg\min}_{Y \in \mathbb{C}^{mp\times p}} ||R_{k} - \mathcal{V}_{m+1}\mathcal{H}_{m}Y||_{F}.
\]
Then the minimal residual solution of the inner method, $A^{-1}A_{\mathcal{C}_{k}}\mathcal{V}_{m}Y$, gives the outer approximation
\[
X_{k+1} =X_{k}+A^{-1}A_{\mathcal{C}_{k}}\mathcal{V}_{m}Y_{m} = X_{k}+(I-\mathcal{U}_{k}\mathcal{C}_{k}^{H}A)\mathcal{V}_{m}Y_{m},
\]
which is also the solution to the global minimization problem
\[
 \min \{||B-AX||_{F}: X\in \mathcal{R}(\mathcal{U}_{k})\oplus \mathcal{R}(\mathcal{V}_{m})\}.
\]
\end{theorem}
\textbf{Proof.}
It is analogous to the proof of Theorem 2.1 \cite{Sturler} and see also \cite{M. Eiermann}. $\Box$

\section{Conclusions and future works}

In this paper, we have derived a new deflated block GCROT($m,k$) method
for nonsymmetric linear systems with multiple RHSs.
Incorporating this modified block Arnoldi deflation, the new algorithm can detect the possible linear dependence at each iteration during the block Arnoldi procedure and avoids expensive computational operations. Moreover, we analyze its main mathematical properties and prove that the deflation procedure is based on a non-increasing behavior of the singular values of the true block residual.
Numerical examples report that the DBGCROT($m,k$) approach can lead to
a faster convergence and is more effective than some other block
solvers, especially when the RHSs are nearly linearly dependent. Therefore, it may be concluded that DBGCROT($m,k$) is a
competitive method for solving the linear systems with multiple
RHSs.\\

 \begin{acknowledgements}
The authors sincerely thank Jason. E. Hicken for suggesting and providing the information of flexible variant of GCROT($m,k$). We are also grateful to the anonymous referees for their valuable and helpful comments that greatly improved the original manuscript of this paper.
 \end{acknowledgements}


\end{document}